\documentclass[final]{siamltex}
\usepackage{amssymb}
\usepackage{amsmath}
\usepackage{lipsum}
\usepackage{subfigure}
\usepackage{ntheorem}
\theorembodyfont{\rmfamily}
\newtheorem{rem}{Remark}
\usepackage{graphicx,color}
\numberwithin{equation}{section}
\usepackage{geometry}

\begin{document}
\title{Asymptotic analysis on the sharp interface limit of the time-fractional Cahn--Hilliard equation}

\author{Tao Tang\thanks{Division of Science and Technology, BNU-HKBU United International College, Zhuhai, Guangdong, China, \&
SUSTech International Center for Mathematics, Southern University of Science and Technology, Shenzhen 518055, China, {\tt tangt@sustech.edu.cn}}
\and
Boyi Wang\thanks{Department of Mathematics, Southern University of Science and Technology, Shenzhen 518055, China, \& Department of Mathematics, National University of Singapore, Singapore 119076, {\tt 11755001@mail.sustech.edu.cn}}
\and
Jiang Yang\thanks{Department of Mathematics \& SUSTech International Center for Mathematics, Southern University of Science and Technology, Shenzhen 518055, China, {\tt yangj7@sustech.edu.cn}}
}

\maketitle
\begin{abstract}
 In this paper, we aim to study the motions of interfaces and coarsening rates governed by the time-fractional Cahn--Hilliard equation (TFCHE). It is observed by many numerical experiments that the microstructure evolution described by the TFCHE displays quite different dynamical processes comparing with the classical Cahn--Hilliard equation, in particular, regarding motions of interfaces and coarsening rates.
 By using the method of matched asymptotic expansions, we first derive the sharp interface limit models. Then we can theoretically analyze the motions of interfaces with respect to different timescales. For instance, for the TFCHE with the constant diffusion mobility, the sharp interface limit model is a fractional Stefan problem at the time scale $t=O(1)$. However, on the time scale $t=O(\varepsilon^\frac1\alpha)$ the sharp interface limit model is a fractional Mullins--Sekerka model. Similar asymptotic regime results are also obtained for the case with one-sided degenerated mobility. Moreover, scaling invariant property of the sharp interface models suggests that the TFCHE with constant mobility preserves an $\alpha/3$ coarsening rate and a crossover of the coarsening rates from $\frac{\alpha}{3}$ to $\frac\alpha4$ is obtained for the case with one-sided degenerated mobility, which are in good agreement with the numerical experiments.
\end{abstract}

\begin{keywords}
{Method of matched asymptotic expansions, time-fractional Cahn--Hilliard equation, phase-field modeling, coarsening rates, motion of interfaces}
\end{keywords}

\begin{AMS}
65M30, 65M15, 65M12
\end{AMS}

\section{Introduction}
The coarsening progress (see Fig.\ref{fig1}-\ref{fig3}) is an ubiquitous phenomena and is observed in many fields such as the study of solid or fluid in material science, opinon dynamic in social science, and pattern formation in biological science \cite{cugliandolo2015coarsening}. It is marked by an increase of the typical length scale in the spatial structures, which is due to the decrease of the interfacial energy \cite{cai2016effect,dai2012motion,dai2016computational,dai2010crossover,dai2005universal,kohn2004coarsening}.
During the coarsening process, a power law, i.e., the increasing of a characteristic length scale with respect to the power of time, is often observed \cite{dai2012motion, dai2016computational, tang2019energy}, as well as Fig.\ref{fig4}. To measure the coarsening process, a coarsening rate is introduced.
It is clear that the Cahn--Hilliard equation (CHE) can be used for simulating the coarsening progress with an 1/3 power law. This power law coincides with the coarsening rate indicated by the classical LSW theory for bulk diffusion.
However, different coarsening rates have also been discovered, which suggests that the CHE is insufficient.
For example, significantly small coarsening rates, i.e. 0.13, 0.07 and 0.09, are observed in the coarsening of $\gamma'$ precipitates \cite{gagaray2017comparative,sequeira1995bimodal}.
The authors explain that as a result of the existence of the elastic strain field.
Besides, a coarsening rate of 1/2 is observed in the study of precipitate in rapidly solidified Al-Si alloy and it is due to a change of the anealling temperature according to the author \cite{cai2016effect}. 
More examples of different coarsening process are introduced in \cite{cugliandolo2015coarsening}. These results suggest that the CHE may not be a suitable coarsening model of every cases and other models should be considered.

Recently, time-fractional models have drawn people's attention \cite{allen2016parabolic,chen2019accurate, du2017analysis, du2020time,jiang2015fast,lin2007finite,liu2018time,le2016numerical,zhang2020non,li2017space}. 
Numerical results have shown that the coarsening rate of a time-fractional Cahn--Hilliard equation (TFCHE) depends not only on the mobility, but also on the order of the fractional derivative \cite{li2017space,zhao2018Time,podlubny1998fractional,tang2019energy,zhao2019power,du2020time}. 
Especially, an intriguing coarsening rate of $\alpha/3$ is observed in \cite{tang2019energy}. Fig. \ref{fig1}-\ref{fig4} show the case for $\alpha = 0.9$. 
%

This paper is concerned with the motion of interfaces and coarsening dynamics  of the time-fractional Cahn--Hilliard equation (TFCHE) 
{\begin{equation}
\left\{
\begin{array}{lr}
\partial_t^\alpha u= \nabla(M(u)\nabla\mu),\\
\mu = -\varepsilon^2\Delta u +F'(u),\quad x\in \Omega,\quad 0<t<T, \label{eq1}
 \end{array}
 \right.
\end{equation}
where, for some given $0<\alpha<1$, $\partial_t^\alpha$ is the Caputo fractional derivative \cite{allen2016parabolic,kilbas2006theory,samko1993fractional} defined by
\begin{displaymath}
\partial_t^\alpha u = \frac{1}{\Gamma(1-\alpha)}\int_0^t\frac{u'(\tau)}{(t-\tau)^\alpha}d\tau,\quad t >0.
\end{displaymath}
As an nonlocal-in-time extension of classical phase-field models, $u$ is the order parameter, $\varepsilon$ represents the width of interfaces, and $\mu$ is the chemical potential. Without loss of generality, we restrict our attention on the commonly used double well potential
\begin{equation}
F(u) = \frac{1}{4}(u^2-1)^2.
\end{equation}
In \eqref{eq1}, the diffusion mobility function $M(u)$ is taken as the constant $1$ or the one-sided degenerate function $1+u$. For simplicity, \eqref{eq1} is subject to the Nuemann boundary conditions
\begin{align}\label{eq2}
\frac{\partial u}{\partial n} = \frac{\partial\mu}{\partial n} =0,\quad x\in\partial\Omega,\quad 0<t<T,
\end{align} and the initial data
\begin{equation}u(x,0) = u_0(x),\quad x\in\Omega.\label{eq3}
\end{equation}


\begin{figure}[h] \centering    
    \subfigure[u at t = 4] {
     \label{fig1}     
    \includegraphics[width=0.48\columnwidth,height=6.3cm]{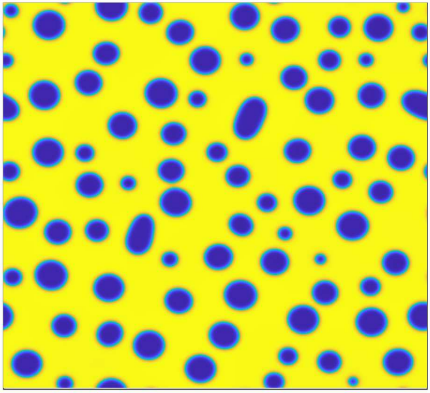}
    }     
    \subfigure[u at t = 25] { 
    \label{fig2}     
    \includegraphics[width=0.48\columnwidth,height=6.3cm]{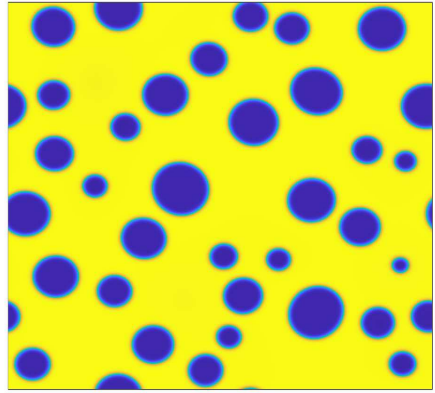}    
    }    
    \subfigure[u at t = 100] { 
    \label{fig3}     
    \includegraphics[width=0.48\columnwidth,height=6.3cm]{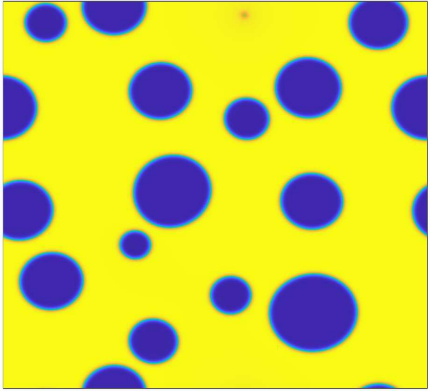}     
    } 
    \subfigure[evolution of energy]{ 
    \label{fig4}     
    \includegraphics[width=0.48\columnwidth,height=6.3cm]{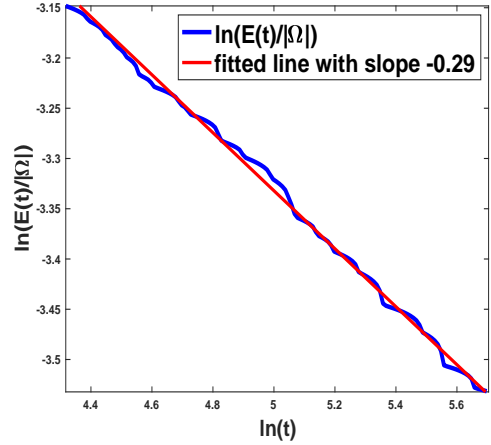}     
    }     
    \caption{$M(u) = 1$, $\alpha = 0.9$, $\varepsilon=0.05$. Morphological patterns at $t=4$ (top left), $t=25$ (top right), $t=100$ (bottom left), $\ln (E(t)/|\Omega|)$ vs. $\ln(t)$ (bottom right).}     
    \label{fig11}     
\end{figure}

Extensive investigations have been made to study the coarsening process and the coarsening rates of the Cahn--Hilliard equations. Pego \cite{pego1989front} studied the asymptotic regimes on CHE with the constant mobility by the method of matched asymptotic expansions. Alikakos, Bates and Chen \cite{alikakos1994convergence} proved the convergence of CHE to the Mullins--Sekerka equations. Cahn, Elliott and Novick-Cohen \cite{cahn1996cahn} studied the degenerate CHE and obtained the surface diffusion model. In addition, it has been shown that coarsening rate of the Cahn--Hilliard equation is related to the diffusion mobility. Dai and Du \cite{dai2012motion,dai2016computational} studied the motion of interfaces for CHE with single-sided degenerate mobility, and they obtained its sharp interface limits as well as the coarsening rates. More results related to the CHE can be found in, i.e., \cite{alikakos1994convergence,bernis1990higher, caffarelli1995bound,dai2010crossover,du2016nonlocal, Jie2010Numerical,jingxue1992existence, dai2005universal} etc.


Motivated by the above asymptotic analysis theory and numerical results on the coarsening rates for time fractional CHE, we will establish asymptotic regime theory on the TFCHE by the method of matched asymptotic expansion as used in \cite{pego1989front} and to derive the surface diffusion models of interface motion for the TFCHE. As far as we know, this is the first work to study the coarsening process and coarsening rate of TFCHEs using formal asymptotic matching.

Our main results are twofold. Firstly, a formal asymptotic description of the TFCHE in the later regime of phase seperation is given, where different types of mobilities are discussed. 
Secondly, using the resulted sharp interface models and the scaling invariant property, we explain the corresponding coarsening rates for the TFCHEs, which agrees well with numerical observations in \cite{tang2019energy,zhao2019power}. 
A more precise outline of the first result is given below. In a slow time scale $O(1)$, the solution at leading order satisfies a nonlocal ``Stefan problem'' with equilibrium condition at the interface, and the leading order inner solution is the solution to the following problem
\begin{subequations}
\begin{align}
 &F'(U)-\partial_{zz}U = 0,\\
 &U(-\infty)=-1, \quad\quad\quad\quad U(+\infty) = 1, \quad U(0) = 0,
\end{align}
\end{subequations}
 which is the re-scaled $\tanh$ function $U(z)=\tanh(z/\sqrt2)$. Then, on a much more slower timescale $t_1=O( \varepsilon^{\frac{1}{\alpha}}t)$, phase equilibrium holds everywhere and interface motion is governed by $\mu_1$, which is the second term in the asymptotic expansion of the chemical potential $\mu$, obeying to the following nonlocal ``Mullins--Sekerka" model
\begin{subequations}
\begin{align}
&& &\partial_{t_1}^\alpha u_0  = \Delta \mu_1, & &  \text{in}\  \Omega\backslash \Gamma,&\label{cm-1}\\
&& &\mu_1 = \kappa \frac{S}{[U]},& & \text{on}\ \Gamma,&\label{cm-2}\\
&& &\text{I}^{1-\alpha} V= [\partial_\textbf{m}\mu]_-^+/[U],& & \text{on}\ \Gamma.&\label{cm-3}
\end{align}
\end{subequations}
 In \eqref{cm-1}-\eqref{cm-3}, $S$ and $[U]$ are some constants, $u_0$ is the sign function of the distance function $\phi$, $\Gamma$ is the interface, $\kappa = \Delta \phi$ is the mean curvature, $V = \partial_{t}\phi$ is the normal velocity of $\Gamma$ on $x$ with the signed distance $\phi$ from the point $x\in \Omega$ to interface, $\textbf{m}$ is the unit normal vector on $\Gamma$, $\text{I}^{1-\alpha}$ denotes the fractional integral operator, 
$u_0$ is determined by the interface $\Gamma$ and equals to $\pm1$ in $\Omega^{\pm}$, correspondingly. The present results reduce to the classical one of Pego \cite{pego1989front} for local CHE
\begin{align*}
&& &V= [\partial_\textbf{m}\mu_{1}]_-^+/[U],& &\text{on}\ \Gamma &&&
\end{align*}
by letting $\alpha\rightarrow 1$.

As for the case with one-sided degenerate mobility, i.e., $M(u) = 1 + u$, the corresponding sharp interface models in time scales $t_1=O( \varepsilon^{\frac{1}{\alpha}}t)$ and $t_2 =O( \varepsilon^{\frac{2}{\alpha}}t)$ are derived respectively as the following nonlocal Mullins--Sekerka models
: 
\begin{subequations}
\begin{align}
&& &\partial_{t_1}^\alpha u_0=\Delta \mu_1,&  &\text{in} \ \Omega^+,&\label{dgms1}\\
&& &\mu_1 =-\kappa \frac{S}{[U]},& &\text{on}\ \Gamma,& \\
&& &\text{I}^{1-\alpha}V = \partial_\textbf{m}\mu_1^+,&  &\text{on}\ \Gamma&\label{dgms2}
\end{align}
\end{subequations}
and
\begin{subequations}
\begin{align}
  && &\partial_{t_2}^\alpha u_0=\nabla(\mu_1\nabla \mu_1),& &\text{in}\ \Omega^-,& \label{dgms21}\\
  && &\mu_1 = -\kappa \frac{S}{[U]}, & &\text{on}\ \Gamma,& \\
 && &2\Delta \mu_2 = \partial_t^\alpha u_0, & &\text{in} \ \Omega^+,& \\
 && &\mu_2 = -\kappa^2\frac{S_1}{[U]},& &\text{on}\ \Gamma,& \\
 && &\text{I}^{1-\alpha}V =\partial_\textbf{m}\mu_2^++\frac14\mu_1^-\partial_\textbf{m}\mu_1^-,& &\text{on}\ \Gamma.& \label{dgms22}
\end{align}
\end{subequations}

A more precise outline of the second model is given below. For the case with the constant mobility $M(u)=1$, 
 the scaling invariant of nonlocal Mullins-Sekerka model implies an coarsening rate of $\alpha/3$, which coincides well with that observed in numerical experiments. For the case with one-sided degenerate mobility $M(u) = 1+u$, the models \eqref{dgms1}-\eqref{dgms2} and \eqref{dgms21}-\eqref{dgms22} exhibit two different coarsening rates of $\frac{\alpha}{3}$ and $\frac{\alpha}{4}$ respectively, which are in good agreement with the observations in \cite{dai2012motion,dai2016computational}.

The rest of the paper is organized as follows. In Section \ref{sec2} and Section \ref{sec3}, we establish sharp interface limit models for the TFCH system \eqref{eq1}-\eqref{eq3} when $M(u) = 1$ and $M(u) = 1+u$, respectively. In Section \ref{sec4}, the scaling invariant properties of sharp interface models and the coarsening rates will be discussed. Some concluding remarks are given in the finial section.

\section{Sharp interface models when $M(u) = 1$} \label{sec2}
The method of matched asymptotic expansions expansion as in Pego \cite{pego1989front} will be used in this section. For all $\gamma \in \mathbb{R}$ and $t_1 = \varepsilon^\gamma t$, simple calculation yields
\begin{align}
\partial_t^\alpha v(t_1) & = \frac{\varepsilon^{\alpha\gamma}}{\Gamma(1-\alpha)}\int_0^{\varepsilon^\gamma t}\frac{v'(\tau)}{(\varepsilon^\gamma t- \tau)^\alpha}d\tau=\varepsilon^{\alpha\gamma}\partial_{t_1}^\alpha v(t_1)\label{9}.
\end{align}
Below we will develop sharp interface models at different time scales. We assume that with domain $\Omega \subset \mathbb{R}^N$, $N=2$ or $3$, there is a smooth $N-1$ dimension interface $\Gamma$ which divides $\Omega$ into $\Omega^+$ and $\Omega^-$, and the interface $\Gamma$ does not intersect with the boundary.

\subsection{The time scale $t = O(1)$: a time-fractional Stefan problem}
We assume that the phase structures are nearly equilibrated.

\subsubsection{Outer expansion} We expand the solution in a series of powers of $\varepsilon$ in the timescale $t$:
\begin{subequations}
\begin{align}
u(x,t) =  u_0(x,t)+\varepsilon  u_1(x,t)+\cdots,\label{eq22a}\\
\mu(x,t) = \mu_0(x,t)+\varepsilon \mu_1(x,t)+\cdots.\label{eq22b}
\end{align}
\end{subequations}
In this time scale,
\begin{align}
    \partial^\alpha_t u =  \partial^\alpha_t u_0+\varepsilon \partial^\alpha_t u_1+\cdots.\label{eq23}
\end{align}
 By comparing Eq. \eqref{eq23} with (\ref{eq1}) and matching the powers of $\varepsilon$, we get
\begin{align}
\partial_{t}^\alpha u_0 = \Delta\mu_0,\quad \mu_0 = F'(u_0).\label{t-eq-1}
\end{align}
This method will be used many times in this paper. The leading order equation implies that the phase parameters evolve according to the chemical potential. The boundary condition on $\partial \Omega$ is taken naturally as $\frac{\partial u_0}{\partial n}=0$. To model this problem, it is necessary to derive the boundary conditions on the interface, which can be done by matching outer solutions with the inner solutions.

\subsubsection{Inner expansion} Now we consider the inner expansions near the front. Intuitively, the inner solutions takes the value of the solutions restricted on the interface. The inner solutions will be defined in this region by an inner variable $z$. Moreover, the inner solution matches with the outer solution when $z\rightarrow \pm\infty$ according to some specified matching conditions. We take the same notations as Pego \cite{pego1989front}. In order to define the inner variable $z$, define the stretched normal distance to the front
\begin{displaymath}
z = \phi(x,t)/\varepsilon,
\end{displaymath}
where $\phi(x,t)$ is the signed distance of the point $x$ in $\Omega$ to the interface $\Gamma(t)$ such that $\phi>0$ in $\Omega^+$ and $\phi<0$ in $\Omega^-$. Note that $\phi$ is a smooth function near $\Gamma$ if $\Gamma$ is smooth.

Consider the functions $\tilde v = \tilde v(z,x,t)$ defined near the interface. Following \cite{pego1989front}, we require that $v$ does not varies when $x$ varies normally to $\Gamma$ but $z$ holds, that is, $\tilde v(z,x+\alpha\nabla\phi,t) =\tilde v(z,x+\nabla\phi,t)$ for small $\alpha$ or $\nabla \phi\cdot\nabla_x \tilde v =0$. Moreover, define
\begin{align}
\textit{\textbf{m}} = \nabla \phi(x,t),\quad \kappa = \Delta \phi(x,t),\quad V(x,t)=\partial_{t}\phi(x,t),
\end{align}
where $\textit{\textbf{m}}$ is the unit normal vector on $\Gamma$ pointing towards $\Omega^+$, $\kappa$ is the mean curvature of $\Gamma$ at point $x$, $\partial_{t}\phi = V(x,t)$ is the normal velocity of front motion in this time scale which is positive when pointing towards $\Omega^-$. We also assume that $\partial_{t}\phi = V(x,t)$ exists  for all $x \in \Omega$. Given $\tilde v(z,x,t)$ and $v = \tilde v(\phi(x,t)/\varepsilon,x,t)$, we have derivatives transform according to the relations \cite{pego1989front}:
\begin{subequations}
\begin{align}
\nabla v &= \nabla_x \tilde v + \varepsilon^{-1}\textbf{m}\partial_z \tilde v,\label{delta1}\\
\Delta v &= \Delta_x \tilde v + \varepsilon^{-1}\kappa\partial_z\tilde v+\varepsilon^{-2}\partial_{zz}\tilde v\label{delta},\\
\partial_t v &= \varepsilon^{-1}\partial_{t}\phi\partial_z \tilde v+\partial_{t} \tilde v\label{delta4}.
\end{align}
\end{subequations}
For the inner expansion, we have
\begin{subequations}
\begin{align}
u(x,t) = \tilde u_0(z,x,t)+\varepsilon \tilde u_1(z,x,t)+\cdots,\label{eq26a}\\
\mu(x,t) = \tilde\mu_0(z,x,t)+\varepsilon \tilde\mu_1(z,x,t)+\cdots.\label{eq26b}
\end{align}
\end{subequations}
By Taylor expansion and \eqref{eq26a}-\eqref{eq26b}, the expansions are related by
\begin{subequations}
\begin{align}
&\tilde\mu_0 = F'(\tilde u_0) -\partial_{zz}\tilde{u_0},\label{eq210}\\
&\tilde\mu_1 = F''(\tilde u_0)\tilde u_1-\partial_{zz} \tilde u_1-\kappa \partial_z \tilde u_0,\\
&\tilde\mu_2 = F''(\tilde u_0)\tilde u_2- \partial_{zz}\tilde u_2 -\kappa \partial z \tilde u_1+\frac12F'''(\tilde u_0)\tilde u_1^2-\Delta_x \tilde u_0.
\end{align}
\end{subequations}
Substituting the expansion back to (\ref{eq1}), using the derivative transform formulas \eqref{delta1}-\eqref{delta4} and matching the lowest order term of $\varepsilon$ shows:
\begin{align}
\partial_{zz}\tilde\mu_0& = 0,\label{20}
\end{align}
Integrating (\ref{20}) and combing \eqref{eq210} derive 
\begin{align}
\tilde\mu_0 = a_0(x,t)z+b_0(x,t)= F'(\tilde u_0) -\partial_{zz}\tilde{u}_0.
\end{align}
Since $\tilde u_0$ must be bounded, $a_0(x,t)$ has to be zero.
Then derive $b_0$ by solving the following system:
\begin{subequations}
\begin{align}
&& &F'(\tilde u_0) -\partial_{zz}\tilde{u}_0= b_0, &\label{eq-xx1} \\
&& &\tilde u_0(+\infty,x,t) = u^+(x,t),\quad \tilde u_0(-\infty,x,t) = u^-(x,t)&\label{xx1}.
\end{align}
\end{subequations}
Letting $z\rightarrow \pm\infty$ in \eqref{eq-xx1} and integrating Eq.\eqref{eq-xx1} with respect to $u$ yield
\begin{subequations}
\begin{align}
&& &F'(u^+(x,t)) = F'(u^-(x,t)) =  b_0(x,t),&\label{eq-2.0}\\
&& &b_0(x,t)(u^+(x,t)-u^-(x,t)) = F(u^+(x,t)) - F(u^-(x,t)).& \label{eq-2.1}
\end{align}
\end{subequations}
Assuming that the leading order inner solution $u_0$ links the two pure phases $\pm 1$, which means
\begin{align}
u^+(x,t) =1,\quad u^-(x,t) =  -1.
\end{align}
Therefore, 
\begin{align}
 b_0 = \tilde\mu_0(z) = 0.
\end{align}
Recall \eqref{eq-xx1} with $b_0=0$. As in \cite{dai2012motion}, we choose the well-known solution profile
\begin{align}
 \tilde u_0(z) =\tanh\left(\frac{z}{\sqrt2}\right) =: U(z). \label{eqU}
\end{align}
Matching with the outer solution by \eqref{t-eq-1} derives the boundary conditions for the equilibrium state
\begin{align}
\mu_0 = 0 \quad \text{on}\quad \Gamma.
\end{align}

For the matching between higher order terms, we follow the ideas provided by Caginalp and Fife in \cite{G1988Dynamics}. Fixing $x$ on $\Gamma$, we seek to match the expansions by requiring formally that
\begin{align}
(\mu_0+\varepsilon\mu_1+\cdots)|_{(x+\varepsilon z\textbf{\textit{m}},t)} \approx (\tilde\mu_0+\varepsilon\tilde\mu_1+\cdots)|_{(z,x,t)},
\end{align}
when $\varepsilon z$ is between $o(1)$ and $O(\varepsilon)$. Expanding the left hand side in powers of $\varepsilon$ as $\varepsilon z\rightarrow 0+$, gives
\begin{align}
\mu_0^+ +\varepsilon(\mu_1^++z{D_m}\mu_0^+) + \varepsilon^2(\mu_2^++z{D_m}\mu_1^+ +\frac{1}{2}z^2D_m^2\mu_0^+)+\cdots,\label{exp-mu}
\end{align}
where $D_m$ denotes the directional derivative along \textbf{\textit{m}} and $\mu_i^+$ is the limit when $z\rightarrow0$ along \textbf{\textit{m}}:
\begin{align}
\mu_i^\pm = \lim_{z\rightarrow 0^\pm} \mu_i(x+z\textbf{\textit{m}},t_1).
\end{align}
Similar results hold for $\varepsilon z\rightarrow0^-$. To match these expansions in \eqref{exp-mu} with the inner expansion, one requires
\begin{subequations}
\begin{align}
&& &\mu_0^\pm(x,t) = \tilde \mu_0(z,x,t),&\quad z\rightarrow \pm\infty,\label{m11}\\
&& &(\mu_1^\pm+z{D_m}\mu_0^\pm)(x,t)  = \tilde \mu_1(z,x,t),&\quad z\rightarrow \pm\infty\label{m21},\\
&& &(\mu_2^++z{D_m}\mu_1^+ +\frac{1}{2}z^2D_m^2\mu_0^+)(x,t)  = \tilde \mu_2(z,x,t),& \quad z\rightarrow \pm\infty.\label{m31}
\end{align}
\end{subequations}

The time derivative in the local frame equals
\begin{align}
&\partial_t^\alpha u(x,t)= \partial_t^\alpha\left((\tilde u_0+\varepsilon \tilde u_1+\cdots)|_{(\phi(x,t)/\varepsilon,x,t)}\right)\notag\\
 =& \frac{1}{\Gamma(1-\alpha)}\int_0^{t}\frac{\varepsilon^{-1}\partial_{\tau}\phi(x,\tau)\partial_z \tilde u_0(\phi(x,\tau)/\varepsilon)}{(t-\tau)^\alpha}d\tau \notag\\
 &+\frac{1}{\Gamma(1-\alpha)}\int_0^{t}\frac{\partial_{\tau}\phi(x,\tau)\partial_z \tilde u_1(\phi(x,\tau)/\varepsilon,x,\tau)+\varepsilon(\partial_{\tau} \tilde u_1(z,x,\tau)|_{z = \phi(x,\tau)/\varepsilon})}{(t-\tau)^\alpha}d\tau+\cdots\notag\\
 =& \frac{1}{\Gamma(1-\alpha)}\int_0^{t}\frac{\varepsilon^{-1}\partial_{\tau}\phi(x,\tau)\partial_z \tilde u_0(\phi(x,\tau)/\varepsilon,x,\tau)}{(t-\tau)^\alpha}d\tau
+h.o.t.
\end{align}
Then matching the $O(\frac1\varepsilon)$ term gives
\begin{align}
\frac{1}{\Gamma(1-\alpha)}\int_0^{t}\frac{\partial_{\tau}\phi(x,\tau)\partial_z \tilde u_0(\phi(x,\tau)/\varepsilon)}{(t-\tau)^\alpha}d\tau  &=\tilde\mu_{1zz}(z,x,t)\label{eqvelo}.
\end{align}
Integrating Eq. \eqref{eqvelo} with respect to $z$ over $(-\infty,+\infty)$, we get
\begin{align}
\frac{1}{\Gamma(1-\alpha)}\int_0^{t}\frac{\phi_\tau(x,\tau)}{(t-\tau)^\alpha}d\tau U\big|_{-\infty}^{+\infty}
 &=\tilde\mu_{1z}\big|_{-\infty}^{+\infty}.
\end{align}
By using the matching condition (\ref{m21}), we derive
\begin{align}
\frac{1}{\Gamma(1-\alpha)}\int_0^{t}\frac{\phi_\tau(x,\tau)}{(t-\tau)^\alpha}d\tau
 &=[\textbf{m}\cdot \nabla\mu_0]_-^+[U]^{-1},\label{eqvelo2}
\end{align}
where $[U] =U\big|_{-\infty}^{+\infty}=2$ and $[\textbf{m}\cdot \nabla\mu_0]_-^+$ denotes the jump of the direction derivative of $\mu$ over the interface along the normal vector. We rewrite Eq. (\ref{eqvelo2}) in the following form using the notation of fractional integral:
\begin{align}
\text{I}^{1-\alpha} V=\frac12[\partial_\textbf{m} \mu_0]_-^+\label{35}.
\end{align}

\textbf{Sharp interface model in $t = O(1)$.} Ignoring the subscripts, the sharp interface model is a time-fractional Stefan model.
\begin{subequations}
\begin{align}
&& &\partial_{t}^\alpha u_0  = \Delta\mu_0,\ \mu_0 = F'(u_0),& &\text{in}\   \Omega/\Gamma,& \\
&& &u_0 = 1,\ \mathrm{on} \quad\Gamma^+,\quad u_0 = -1,\ \mathrm{on} \quad \Gamma^-,& \\
&& &\mu_0 = 0,& &\text{on}\ \Gamma,&\\
&& &\text{I}^{1-\alpha} V = \frac12[\partial_\textbf{m}\mu_0]_-^+,&
\end{align}
\end{subequations}

\subsection{The time scale $t_1 = \varepsilon^{\frac{1}{\alpha}}t$: a time-fractional Mullins--Sekerka model}
In this part we derive the time-fractional sharp interface model in the time scale $t_1 = \varepsilon^{\frac{1}{\alpha}}t$.\\

\subsubsection{Outer expansion} In this time scale,
\begin{align}
    \partial^\alpha_t u =  \varepsilon\partial^\alpha_{t_1} u_0+\varepsilon^2 \partial^\alpha_{t_1} u_1+\cdots.\label{eq24}
\end{align}
Similar to \eqref{t-eq-1}, we have
\begin{align}
0 = \Delta\mu_0,\quad \mu_0 = F'(u_0), \quad \partial_{t_1}^\alpha u_0 = \Delta\mu_1,\quad \partial_{t_1}^\alpha u_1 = \Delta\mu_2.\label{sharpeq}
\end{align}
In this time scale, at leading order, we have a steady state equation for $\mu_0$. Nevertheless, Eq.\eqref{sharpeq} the conditions on the interface and the boundary are also required. The boundary condition on $\partial\Omega$ is naturally inherited from the boundary condition $\frac{\partial u}{\partial n}=0$, but for the boundary conditions on the interface, we need to solve for them by asymptotically matching the outer solutions and the inner solutions.\\

\subsubsection{Inner expansion} Similar to \eqref{20}, matching $1/\varepsilon^2$ and $1/\varepsilon$ terms in the second equation in Eq. \eqref{eq1} yields
\begin{subequations}
\begin{align}
&& &\partial_{zz} \tilde \mu_0 = 0,\label{11}\\
&& &\kappa \partial_z \tilde\mu_0 +\partial_{zz}\tilde\mu_1 = 0.\label{22}
\end{align}
\end{subequations}

Analogous analysis to section 2.1 leads to a tanh profile again, i.e.,
\begin{displaymath}
\tilde u_0 = U(z), \quad \tilde\mu_0 = 0,
\end{displaymath}
where $U(z)$ is defined in \eqref{eqU}. We also assume that $u_0^+(x,t_1) = 1,\quad u_0^-(x,t_1) = -1$. Matching the inner solution with the outer solution according to (\ref{m11}), one derives the boundary conditions for the outer solution
\begin{align}
&& &\mu_0 = 0,\quad\text{on}\ \Gamma.&
\end{align}
Notice that now $\Delta \mu_0 = 0$ and $\mu_0 = F'(u_0)$ in \eqref{sharpeq}, therefore,
\begin{displaymath}
\mu_0 = 0 \quad \text{in} \quad\Omega,\quad  u_0 \equiv -1 \quad \text{in}\quad \Omega^-;\quad \text{and}\quad u_0 \equiv 1 \quad \text{in} \quad \Omega^+.
\end{displaymath}
As for $\tilde\mu_1$, we have
\begin{equation}
\tilde\mu_1 = F''(\tilde u_0)\tilde u_1-\partial_{zz} \tilde u_1-\kappa \partial_z \tilde u_0=b_2(x,t_1).\label{tt}
\end{equation}
Since $ F''(\tilde u_0)\tilde u_0'(z) -\partial_{zz}\tilde u'_0= 0$,  multiplying (\ref{tt}) by $U'$ and integrating by $z$ on $(-\infty,+\infty)$ yield
\begin{displaymath}
[U]\tilde\mu_1 + \kappa S = 0,
\end{displaymath}
where $$\quad S = \int_{-\infty}^{+\infty}U'(z)^2dz, \quad [U] = u^+-u^-=2.$$
Using the matching conditions (\ref{m21}),
\begin{align*}
&& &\mu_1 = \tilde \mu_1 =-\kappa \frac{S}{[U]},\quad\text{on}\ \Gamma.&
\end{align*}
Letting $\varepsilon \rightarrow 0$, we have the boundary conditions of $\mu$ at the interface $\Gamma$. Therefore, we have a closed system for $\mu_1$
\begin{subequations}
\begin{align}
&& &\partial_{t_1}^\alpha u_0 = \Delta \mu_1,&  &\text{in} \ \Omega\backslash \Gamma,& \label{sharpmullin1}\\
&& &\mu_1 = -\kappa \frac{S}{[U]},& &\text{on}\ \Gamma,& \\
&& &\partial_\textbf{m} \mu_1= 0,& &\text{on}\ \partial\Omega.& \label{sharpmullin3}
\end{align}
\end{subequations} Provided that $\Gamma$ is known and smooth, which is well-defined and can be solved independently in each $\Omega^\pm$.\\

Similar to (\ref{eqvelo2}), in this new time scale we have
\begin{align}
\frac{1}{\Gamma(1-\alpha)}\int_0^{t_1}\frac{\phi_\tau(x,\tau)}{(t_1-\tau)^\alpha}d\tau
 &=[\textbf{m}\cdot \nabla\mu_1]_-^+[U]^{-1},
\end{align}
which is
\begin{align}
\text{I}^{1-\alpha}V = \frac12[\partial_\textbf{m} \mu_1]_-^+.
\end{align}
\textbf{Sharp interface model in $t_1 = \varepsilon^{\frac{1}{\alpha}}t$.} Collecting the above equations (\ref{sharpmullin1})-(\ref{sharpmullin3}), we get the sharp interface model as follows
\begin{subequations}
\begin{align}
&& &\partial_{t_1}^\alpha u_0 = \Delta \mu_1,& &\text{in} \ \Omega\backslash \Gamma,&\label{sharpmullin21}\\
&& &\mu_1 = -\kappa \frac{S}{[U]},& &\text{on}\ \Gamma,&\\
&& &\frac{\partial\mu_1}{\partial n} = 0,& &\text{on} \ \partial\Omega,&\\
&& &\text{I}^{1-\alpha}V = \frac12[\partial_\textbf{m}\mu_1]_-^+,& &\text{on}\ \Gamma.&\label{sharpmullin24}
\end{align}
\end{subequations}
$u_0 \equiv 1$ or $u_0 \equiv 1$ when $\phi>0$ or $\phi<0$, respectively. The system (\ref{sharpmullin21})-(\ref{sharpmullin24}) is well-posed, which determines the motion of the front for given smooth initial data. It is a time-fractional Mullins--Sekerka model.\\

\begin{rem}
Since $u_0$ is the sign function of $\phi$, $\partial_{t} u_0\equiv0$. Hence, in the local CH model, \eqref{sharpmullin21} becomes
$$\Delta \mu_1=0 \quad \text{in} \ \Omega\backslash \Gamma.$$
But for the TFCHE, it is necessary to keep $\partial_{t_0}^\alpha u_0$ due to the nonlocal effect. 
\end{rem}

\section{Sharp interface models when $M(u) = 1+u$}\label{sec3} 

In this section, we intend to derive the sharp interface models of the TFCHE with one-sided mobility $M(u) = 1+u$ under the same problem setting as in Section \ref{sec2}.
To begin with, special treatments are required for the degenerate mobility since in this case the leading order term $1+u_0$ in the asymptotic expansion of $M(u)$ might not be valid when $z\rightarrow-\infty$. Assuming that $1+\tilde u_0$ decreases exponentially, that is $1+\tilde u_0\sim e^{z/\sigma}, z\rightarrow -\infty$. Taking $\eta = \sigma\ln\frac1\varepsilon$, we have the following estimates of $1+\tilde u_0$:
\begin{equation}
1+\tilde u_0=\left\{
\begin{array}{c}
O(\varepsilon),\quad \text{if}\quad z\leq-\eta,\\
O(\varepsilon^2),\quad \text{if}\quad z\leq-2\eta,\\
O(\varepsilon^3),\quad \text{if}\quad z\leq-3\eta,\\
O(\varepsilon^4),\quad \text{if}\quad z\leq-4\eta.
\end{array}
\right.
\end{equation}
To simplify the notations, we denote $\chi_4 = 1_{(-\infty,-4\eta]}$, $\chi_3 = 1_{(-4\eta,-3\eta]}$, $\chi_2 = 1_{(-3\eta,2\eta]}$, $\chi_1 = 1_{(-2\eta,\eta]}$ and $\chi_0 = 1_{(-\eta,+\infty)}$, which are the corresponding characteristic functions on each interval. Then we have the following expansion of $1+\tilde u_0$:
\begin{align}
&1+\tilde u_0 \notag\\
= &(1+\tilde u_0)\chi_0+\varepsilon(1+\tilde u_0)\varepsilon^{-1}\chi_1+\varepsilon^2(1+\tilde u_0)\varepsilon^{-2}\chi_2+\varepsilon^3(1+\tilde u_0)\varepsilon^{-3}\chi_3+\varepsilon^4(1+\tilde u_0)\varepsilon^{-4}\chi_4.
\end{align}
Replacing $1+\tilde u_0$ by the above expansion gives a valid series of $M(u)$. Moreover, similar idea is applied for $\tilde u_{0z}$. When $z\rightarrow -\infty$, $\tilde u_{0z}$ decays at the same rate as $1+\tilde u_0$. As for when $z\rightarrow+\infty$, we assume that $\tilde u_{0z}\sim e^{-z/\hat\sigma}$ and $\hat\eta = \hat\sigma\ln\frac1\varepsilon$ so $\tilde u_{0z}\leq O(\varepsilon)$ when $z\geq\hat\eta$.
Let the partitions be $[-\eta,\hat\eta)$, $[-2\eta,-\eta)\cup[\hat\eta,2\hat\eta)$, $[-3\eta,-2\eta)\cup[2\hat\eta,3\hat\eta)$,$[-4\eta,-3\eta)\cup[3\hat\eta,4\hat\eta)$ and $(-\infty,-4\eta)\cup[4\hat\eta,+\infty)$ and the corresponding characteristic functions be $\hat\chi_0, \hat\chi_1, \hat\chi_2, \hat\chi_3, \hat\chi_4$. Then the following expansion holds:
\begin{align}
\tilde u_{0z} = \tilde u_{0z}\hat\chi_0+\varepsilon\tilde u_{0z}\varepsilon^{-1}\hat\chi_1+\varepsilon^2\tilde u_{0z}\varepsilon^{-2}\hat\chi_2+\varepsilon^3\tilde u_{0z}\varepsilon^{-3}\hat\chi_3+\varepsilon^4\tilde u_{0z}\varepsilon^{-4}\hat\chi_4.
\end{align}

With above expansions, we can now compute $\nabla\cdot(M(u)\nabla \mu)$ as follows
\begin{align}
\nabla\cdot(M(u)\nabla \mu) = M'(u)\nabla_x u\cdot\nabla_x \mu+\varepsilon^{-2}\partial_z u\partial_z \mu+M(u)(\Delta_x \mu+\varepsilon^{-1}\kappa\partial_z\mu+\varepsilon^{-2}\partial_{zz}\mu)\label{mumu}.
\end{align}
By simple calculations, we find the terms of powers of $\varepsilon$ in (\ref{mumu}) correspondingly. The first four leading order terms are required in our later analysis, which are the $O(\frac{1}{\varepsilon^2})$ term
\begin{align}
\hat\chi_0\tilde u_{0z}\tilde \mu_{0z} +\chi_0(1+\tilde u_0)\partial_{zz}\tilde\mu_0,
\end{align}
the $O(\frac{1}{\varepsilon})$ term
\begin{align}
&(\tilde u_{0z}\hat\chi_1\varepsilon^{-1}+\tilde u_{1z})\tilde\mu_{0z}+\tilde u_{0z}\hat\chi_0\tilde \mu_{1z}+\chi_0\kappa(1+\tilde u_0)\tilde \mu_{0z}+\chi_0(1+\tilde u_0)\tilde \mu_{1zz}\\
+&(\tilde u_1+\chi_1(1+\tilde u_0)\varepsilon^{-1})\tilde \mu_{0zz}\notag,
\end{align}
the $O(1)$  term
\begin{align}
&\nabla_x\tilde u_0 \nabla_x\tilde \mu_0+(\hat \chi_0\tilde u_{0z}\tilde \mu_{2z}+(\hat \chi_1\tilde u_{0z}\varepsilon^{-1}+\tilde u_{1z})\tilde \mu_{1z}+(\hat \chi_2\tilde u_{0z}\varepsilon^{-2}+\tilde u_{2z})\tilde\mu_{0z})\\
+&\chi_0(1+\tilde u_0)(\Delta_x\tilde\mu_0+\kappa\tilde\mu_{1z}+\tilde\mu_{2zz})+(\varepsilon^{-1}(1+\tilde u_0)\chi_1+\tilde u_1)(\kappa\tilde\mu_{0z}+\tilde\mu_{1zz})\notag\\
+&(\varepsilon^{-2}(1+\tilde u_0)\chi_2+\tilde u_2)\tilde\mu_{0zz}\notag,
\end{align}
and the $O(\varepsilon)$  term
\begin{align}
&\nabla_x\tilde u_0\nabla_x\tilde \mu_1+\nabla_x\tilde u_1\nabla_x\tilde\mu_0+(\hat\chi_0\tilde u_{0z}\tilde \mu_{3z}+(\hat \chi_1\tilde u_{0z}\varepsilon^{-1}+\tilde u_{1z})\tilde\mu_{2z}\\
+&(\hat \chi_2\tilde u_{0z}\varepsilon^{-2}+\tilde u_{2z})\tilde\mu_{1z}+(\hat \chi_3\tilde u_{0z}\varepsilon^{-3}+\tilde u_{3z})\tilde\mu_{0z})+\notag\\
+&\chi_0(1+\tilde u_0)(\Delta_x\tilde\mu_1+\kappa\tilde\mu_{2z}+\tilde\mu_{3zz})\notag\\
+&((1+\tilde u_0)\chi_1\varepsilon^{-1}+\tilde u_1)(\Delta_x\tilde\mu_0+\kappa\tilde\mu_{0z}+\tilde \mu_{2zz})\notag\\
+&((1+\tilde u_0)\chi_2\varepsilon^{-2}+\tilde u_2)(\kappa\tilde \mu_{0z}+\tilde \mu_{1zz})+((1+\tilde u_0)\chi_3\varepsilon^{-3}+\tilde u_3)\tilde\mu_{0zz}\notag.
\end{align}

We start with a non-trivial time scale in this section.
\subsection{The time scale $t = O(1)$: a one-sided time-fractional Stefan problem}
\subsubsection{Outer expansion}
Similar to \eqref{sharpeq}, it yields
\begin{align}
\partial_t^\alpha u_0  = \nabla((1+u_0)\nabla\mu_0),\quad\partial_t^\alpha u_1 = \nabla((1+u_0)\nabla\mu_1+u_1\nabla\mu_0).\label{eq-u0-11}
\end{align}

\subsubsection{Inner expansion}
In the same way as \eqref{11}- \eqref{22}, the $O(\varepsilon^{-2})$ equation is
\begin{align}
0&=\hat\chi_0\tilde u_{0z}\tilde \mu_{0z} +\chi_0(1+\tilde u_0)\tilde\mu_{0zz}.\label{eq312}
\end{align}
We rewrite Eq.\eqref{eq312} in the following form
\begin{equation}
\chi_0\partial_{z}((1+\tilde u_0)\tilde\mu_{0z})+\hat\chi_0(1+\tilde u_0)\tilde\mu_{0zz}=0.
\end{equation}
That is, for $z\in(-\eta,\hat\eta)$,
\begin{align}
\tilde\mu_{0zz}(1+\tilde u_0)+\tilde \mu_{0z}\tilde u_{0z}=\partial_z(\tilde\mu_{0z}(1+\tilde u_0))=0,
\end{align}
which implies $\tilde \mu_{0z}(1+\tilde u_0) = c_1$ in $(-\eta,\hat\eta)$ and $c_1$ is a constant independent of $z$.
For $z$ in $[\hat\eta,+\infty)$, we have
\begin{equation}
\tilde\mu_{0zz}(1+\tilde u_0) = 0.
\end{equation}
In this case, $\tilde\mu_0 = a_1z+b_1$. Here $a_1$ and $ b_1$ are some functions independent of $z$. However, we claim $a_1 =0$ since $\tilde\mu_0$ must be bounded. This leads to $\tilde \mu_0 = b_1$. Moreover, recall that
\begin{equation}\tilde\mu_0 = F'(\tilde u_0)-\partial_{zz}\tilde u_0,
\quad \tilde u_0|_{x = \pm\infty} =\pm1 ,\label{eq-u0-22}\end{equation}
we take the profile
\begin{equation}
\tilde u_0 =\tanh(z/\sqrt2), \quad \tilde \mu_0 = 0,\quad \forall z \in [\hat\eta,+\infty).\label{tan-fun}
\end{equation}
By the smooth continuity of $\tilde\mu_0$ at $\hat\eta$, we have $c_1 = 0$ and $\tilde\mu_0 = 0$ in $(-\eta,+\infty)$.


Now we consider the governing function of the front. The time-fractional derivative in this scaling is
\begin{align}
\partial_t^\alpha u(x,t)=& \partial_t^\alpha\big((\tilde u_0+\varepsilon \tilde u_1+\cdots)(\phi(x,t)/\varepsilon,x,t)\big)\label{frac2}\\
 =& \frac{1}{\Gamma(1-\alpha)}\int_0^{t}\frac{\varepsilon^{-1}\partial_{\tau}\phi(x,\tau)\partial_z \tilde u_0(\phi(x,\tau)/\varepsilon,x,\tau)+(\partial_{\tau} \tilde u_0(z,x,\tau)|_{z = \phi(x,\tau)/\varepsilon})}{(t-\tau)^\alpha}d\tau \notag\\
 &+\frac{1}{\Gamma(1-\alpha)}\int_0^{t}\frac{\partial_{\tau}\phi(x,\tau)\partial_z \tilde u_1(\phi(x,\tau)/\varepsilon,x,\tau)+\varepsilon\partial_{\tau} \tilde u_1(z,x,\tau)|_{z = \phi(x,\tau)/\varepsilon}}{(t-\tau)^\alpha}d\tau+\cdots\notag.
\end{align}
By matching the $O(1/\varepsilon)$ terms in equation (\ref{frac2}) together with (\ref{mumu}), we yields the following equation
\begin{align}
&\frac{1}{\Gamma(1-\alpha)}\int_0^{t}\frac{\partial_{\tau}\phi(x,\tau)\partial_z \tilde u_0(\phi(x,\tau)/\varepsilon)\chi_0}{(t-\tau)^\alpha}d\tau\label{eq-u0}\\
=&(\tilde u_{0z}\hat\chi_1\varepsilon^{-1}+\tilde u_{1z})\tilde\mu_{0z}+\tilde u_{0z}\hat\chi_0\tilde\mu_{1z}+\chi_0\kappa(1+\tilde u_0)\tilde\mu_{0z}\notag\\
&+\chi_0(1+\tilde u_0)\tilde\mu_{1zz}+\notag
(\tilde u_1+\chi_1\tilde u_0\varepsilon^{-1})\tilde\mu_{0zz}\notag\\
=& \chi_0(1+\tilde u_0)\tilde \mu_{1zz}+\hat\chi_0 \tilde u_{0z}\tilde\mu_{1z}\notag\\
=&\chi_0\partial_z((1+\tilde u_0)\tilde \mu_{1z})+(\hat\chi_0-\chi_0) \tilde u_{0z}\tilde\mu_{1z},\notag
\end{align}
which is simplified by using the former results \eqref{tan-fun}.
Integrating the equation \eqref{eq-u0} over $(-\infty,\infty)$, we have
\begin{align}
\partial^\alpha_{t_1}\phi \tilde u_{0}|_{-\eta}^{\hat\eta} =& (1+\tilde u_0)\tilde \mu_{1z}|_{-\eta}^{+\infty}-\tilde u_{0z}\tilde\mu_{1z}|_{\hat\eta}^{+\infty}.
\end{align}
Here $-1 \leq \tilde u_0(-\eta)\leq -1+O(\varepsilon)$ and $1-O(\varepsilon)\leq \tilde u_0(\hat\eta)\leq 1$. In addition, since $\tilde\mu_{1z}^{+\infty}=0$, we could derive
\begin{align}
\partial^\alpha_t\phi(2+O(\varepsilon)) =  2\lim_{z\rightarrow+\infty}\tilde\mu_{1z}- (1+\tilde u_0)\tilde \mu_{1z}|_{-\eta}+\tilde u_{0z}\tilde\mu_{1z}|_{\hat\eta}=2\lim_{z\rightarrow+\infty}\tilde\mu_{1z}+O(\varepsilon).
\end{align}
Therefore, using the matching conditions,
\begin{align}
\partial^\alpha_t\phi(2+O(\varepsilon)) =  2\lim_{z\rightarrow+\infty}\tilde\mu_{1z}+O(\varepsilon) = 2\partial_\textbf{m}\tilde\mu_0^++O(\varepsilon).
\end{align}
By leting $\varepsilon \rightarrow 0$, we derive the sharp interface condition:
\begin{align}
\partial^\alpha_t\phi = \partial_\textbf{m}\tilde\mu_0^+.\label{eq-u0-33}
\end{align}
\textbf{Sharp interface model in $t = O(1)$.} Combining \eqref{eq-u0-11}, \eqref{eq-u0-22} and \eqref{eq-u0-33}, we derive the following sharp interface model
\begin{subequations}
\begin{align}
&& &\partial^\alpha_t u_0= \nabla((1+u_0)\nabla\mu_0),\quad \mu_0 = F'(u_0),& &\text{in}\quad\Omega\slash\Gamma&,\\
&& &u_0= \pm1\quad \text{on}\quad \Gamma^\pm,\quad \mu_0=0 \mathrm{on}\quad\Gamma,&\\
&& &\text{I}^{1-\alpha} V = \partial_\textbf{m}\mu_0^+,& &\text{on}\quad \Gamma.&
\end{align}
\end{subequations}

\subsection{The time scale $t_1 =\varepsilon^{\frac1\alpha} t$: a one sided time-fractional Mullins--Sekerka(MS) model}
\subsubsection{Outer expansion}
The same as \eqref{eq-u0-11}, by asymptotic matching it yields
\begin{subequations}
\begin{align}
&& &0 = \nabla((1+u_0)\nabla\mu_0),&\\
&& &\partial_{t_1}^\alpha u_0 = \nabla((1+u_0)\nabla\mu_1+u_1\nabla\mu_0),&\label{eq-mu-01}\\
&& &\partial_{t_1}^\alpha u_1 = \nabla((1+u_0)\nabla\mu_2+u_1\nabla\mu_1+u_2\nabla\mu_0).&
\end{align}
\end{subequations}
Here $\mu_0, \mu_1, \mu_2$ are the same as shown before. The first equation implies a equilibrium state, so we take the following solution in the outer region
\begin{equation}
u_0 =\left \{
\begin{array}{c}
+1, \quad \text{in} \quad \Omega^+\\
-1, \quad \text{in} \quad \Omega^-.
\end{array}\right.
\end{equation}

\subsubsection{Inner expansion}
The same as \eqref{eq312}, asymptotic matching leads to
\begin{align}
0=&\tilde u_{0z}\tilde \mu_{0z} +\chi_0(1+\tilde u_0)\tilde\mu_{0zz},\label{eq-u0-41}\\
0=&(\tilde u_{0z}\hat\chi_1\varepsilon^{-1}+\tilde u_{1z})\tilde\mu_{0z}+\tilde u_{0z}\hat\chi_0\tilde\mu_{1z}+\chi_0\kappa(1+\tilde u_0)\tilde\mu_{0z}+\chi_0(1+\tilde u_0)\tilde\mu_{1zz}+\notag\\
&(\tilde u_1+\chi_1\tilde u_0\varepsilon^{-1})\tilde\mu_{0zz}.\label{eq-u0-42}
\end{align}
Now we solve the equations \eqref{eq-u0-41}-\eqref{eq-u0-42} as follows. The equation \eqref{eq-u0-41} is
\begin{equation}
\chi_0\partial_{z}((1+\tilde u_0)\tilde\mu_{0z})+\hat\chi_0(1+\tilde u_0)\tilde\mu_{0zz}=0.\label{eqref}
\end{equation}
For $z\in(-\eta,\hat\eta)$, the equation \eqref{eqref} is
\begin{align}
\tilde\mu_{0zz}(1+\tilde u_0)+\tilde \mu_{0z}\tilde u_{0z} &=\partial_z(\tilde\mu_{0z}(1+\tilde u_0))=0,\label{eq-u0-43}
\end{align}
and, for $z\in[\hat\eta,+\infty)$, it is
\begin{align}
\tilde\mu_{0zz}(1+\tilde u_0) = 0.\label{eq-u0-44}
\end{align}
As in section 3.1, the equations \eqref{eq-u0-43}-\eqref{eq-u0-44} can be solved by the exact function
\begin{align}
\tilde u_0 = \tanh(z/\sqrt2),\quad \tilde \mu_0 = 0.\label{ssss}
\end{align}

Next, we intend to determine $\tilde u_1$ and $\tilde \mu_1$. Using \eqref{ssss}, the equation \eqref{eq-u0-42} is simplified into
\begin{align}
\chi_0(1+\tilde u_0)\tilde \mu_{1zz}+\hat\chi_0 \tilde u_{0z}\tilde\mu_{1z}=0,
\end{align}
which implies that $\tilde\mu_1 = c_2$ for $z\in(\hat\eta,+\infty)$. We assume that $\tilde\mu_1=c_2$ for all $z$. Here $c_2$ is some constant independent of $z$.
Recall that
\begin{equation}
\tilde\mu_1 =-\tilde u_{1zz}-\kappa\tilde u_{0z}+F''(\tilde u_0)\tilde u_1.\label{eq-u1-11}
\end{equation}
Noticing that $F''(u_0)u_0'-\partial_{zz}u_0' = 0$, multiplying the equation \eqref{eq-u1-11} by $u_0'$ and integrating the resulting one over $(-\infty,+\infty)$, we have
\begin{equation}
\tilde\mu_1=c_2 = -\kappa \frac{S}{[U]},\label{eq-mu-02}
\end{equation}
where $S= \int_{-\infty}^{+\infty}\tilde u_0'(z)^2dz$ and $[U] = \tilde u_0|^{+\infty}_{-\infty}$.

As for $\tilde u_1$, we use the idea which was presented in \cite{dai2012motion}. We find that $\tilde u_1 = \kappa\Phi_0+\alpha \tilde u_0'$, where $\Phi_0$ satisfies
\begin{align} 
    -\Phi_{0zz}+F''(\tilde u_0)\Phi_0= \tilde u_{0z}-\frac{S}{[U]}.
\end{align}
 We impose $\tilde u_1(0)=0$ to center the function. Thus it is determined that \begin{align}
    \tilde u_1  = \kappa \Phi= \kappa(\Phi_0-\frac{\Phi_0(0)}{\tilde u'_0(0)}\tilde u_0'),
\end{align}
 where $\Phi(\pm\infty) = -\frac{S}{[U]F''(\pm1)}$.

Now we derive the equation of the front line. Matching with respect to series of $\varepsilon$, we get
\begin{align}
&\frac{1}{\Gamma(1-\alpha)}\int_0^{t_1}\frac{\partial_{\tau}\phi(x,\tau)\partial_z \tilde u_0(\phi(x,\tau)/\varepsilon)\chi_0}{(t_1-\tau)^\alpha}d\tau\\
=&\nabla_x\tilde u_0 \nabla_x\tilde \mu_0+((\hat \chi_0\tilde u_{0z}\varepsilon^{-1}+\tilde u_{1z})\tilde\mu_{1z}+\hat \chi_0\tilde u_{0z}\tilde \mu_{2z}+(\hat \chi_2\tilde u_0\varepsilon^{-2}+\tilde u_{2z})\tilde\mu_{0z})\notag\\
&+\chi_0(1+\tilde u_0)(\Delta_x\tilde\mu_0+\kappa\tilde\mu_{1z}+\tilde\mu_{2zz})(\varepsilon^{-1}(1+\tilde u_0)\chi_1+\tilde u_1)(\kappa\tilde\mu_{0z}+\tilde\mu_{1zz})+\notag\\
&+(\varepsilon^{-2}(1+\tilde u_0)\chi_2+\tilde u_2)\tilde\mu_{0zz},\notag
\end{align}
which yields, by using the known functions $\tilde u_0, \tilde \mu_0, \tilde u_1, \tilde \mu_1$, that
\begin{align}
\frac{1}{\Gamma(1-\alpha)}\int_0^{t_1}\frac{\varepsilon^{-1}\partial_{\tau}\phi(x,\tau)\partial_z \tilde u_0(\phi(x,\tau)/\varepsilon)\chi_0}{(t_1-\tau)^\alpha}d\tau
=\hat \chi_0\tilde u_{0z}\tilde\mu_{2z}+\chi_0(1+\tilde u_0)\tilde\mu_{2zz},
\end{align}
which gives, by integrating in$(-\infty,\infty)$ and using the matching condition, that
\begin{align}
\partial^\alpha_{t_1}\phi =  \lim_{z\rightarrow+\infty}\tilde\mu_{2z} = \partial_\textbf{m}\mu_1^+.\label{eq-mu-03}
\end{align}

\textbf{Sharp interface model in $t_1 = \varepsilon^{\frac{1}{\alpha}}t$.} It follows from \eqref{eq-mu-01},\eqref{eq-mu-02} and \eqref{eq-mu-03} that the sharp interface model in this timescale is
\begin{subequations}
\begin{align}
&& &\partial^\alpha_{t_1} u_0 = \Delta \mu_1,& &\text{in} \quad\Omega^+,&\label{eq-mu-001}\\
&& &\mu_1 = -\kappa \frac{S}{2},& &\text{on}\quad \Gamma,&\label{eq-mu-002}\\
&& &\text{I}^{1-\alpha}V = \partial_\textbf{m}\mu_1^+,&  &\text{on}\quad \Gamma&\label{eq-mu-003}.
\end{align}
\end{subequations}
$u_0$ is the sign function of $\phi$ and $u_0 = \pm1$ in $\Omega^\pm$. We call \eqref{eq-mu-001}-\eqref{eq-mu-003} the time-fractional MS model. The front motion is governed only by the phase parameter restricted in $\Omega^+$.

\subsection{The time scale $t_2 =\varepsilon^{\frac2\alpha} t$}
\subsubsection{Outer expansion}
In this case, by asymptotic matching, it yields
\begin{subequations}
\begin{align}
0& = \nabla((1+u_0)\nabla\mu_0),\label{eq-t2-1}\\
0&= \nabla((1+u_0)\nabla\mu_1+u_1\nabla\mu_0),\label{eq-t2-2}\\
\partial_{t_2}^\alpha u_0 &= \nabla((1+u_0)\nabla\mu_2+u_1\nabla\mu_1+u_2\nabla\mu_0).\label{eq-t2-3}
\end{align}
\end{subequations}
Let us solve the equations \eqref{eq-t2-1}-\eqref{eq-t2-3}.
The equation \eqref{eq-t2-1} implies a equilibrium state, so it is reasonable to take static solutions in $\Omega^+$ and $\Omega^-$
\begin{equation}
u_0 =\left \{
\begin{array}{c}
+1, \quad\quad \text{in} \quad \Omega^+,\\
-1, \quad\quad \text{in} \quad \Omega^-,
\end{array}\right.
\end{equation}
which yields, together with \eqref{eq-t2-2}-\eqref{eq-t2-3}, that the governing equations of $\mu_1$ in $\Omega^-$and $\mu_2$ in $\Omega^+$ are
\begin{subequations}
\begin{align}
&& &\nabla(\mu_1\nabla \mu_1) = \partial_{t_2}^\alpha u_0,& \text{in} \quad \Omega^-,\label{outer31}\\
&& &\Delta \mu_1 = 0,& \text{in} \quad \Omega^+,\label{outer311}\\
&& &2\Delta \mu_2 + \frac12\nabla(\mu_1\nabla \mu_1) = \partial_{t_2}^\alpha u_0,& \text{in} \quad \Omega^+.\label{outer32}
\end{align}
\end{subequations}
Therefore, we take the solution that $\mu_1$ is a constant in $\Omega^+$.
\subsubsection{Inner expansion}
Similarly, assymptotic matching $\varepsilon$ yields
\begin{align}
0=&\tilde u_{0z}\tilde \mu_{0z} +\chi_0(1+\tilde u_0)\tilde\mu_{0zz},\\
0=&(\tilde u_{0z}\hat\chi_1\varepsilon^{-1}+\tilde u_{1z})\tilde\mu_{0z}+\tilde u_{0z}\hat\chi_0\mu_{1z}+\chi_0\kappa(1+\tilde u_0)\tilde\tilde\mu_{0z}+\chi_0(1+\tilde u_0)\tilde\mu_{1zz}\\
&+(\tilde u_1+\chi_1\tilde u_0\varepsilon^{-1})\tilde\mu_{0zz},\notag\\
0=&\nabla_x\tilde u_0 \nabla_x\tilde \mu_0+((\hat \chi_0\tilde u_{0z}\varepsilon^{-1}+\tilde u_{1z})\tilde\mu_{1z}+\hat \chi_0\tilde u_{0z}\tilde \mu_{2z}+(\hat \chi_2\tilde u_0\varepsilon^{-2}+\tilde u_{2z})\tilde\mu_{0z})\label{eq-u3-01}\\
&+\chi_0(1+\tilde u_0)(\Delta_x\tilde\mu_0+\kappa\tilde\mu_{1z}+\tilde\mu_{2zz})+(\varepsilon^{-1}(1+\tilde u_0)\chi_1+\tilde u_1)(\kappa\tilde\mu_{0z}+\tilde\mu_{1zz})\notag\\
&+(\varepsilon^{-2}(1+\tilde u_0)\chi_2+\tilde u_2)\tilde\mu_{0zz},\notag
\end{align}
and
\begin{align}
&\frac{1}{\Gamma(1-\alpha)}\int_0^{t_2}\frac{\partial_{\tau}\phi(x,\tau)\partial_z \tilde u_0(\phi(x,\tau)/\varepsilon)\chi_0}{(t_2-\tau)^\alpha}d\tau\label{eqphi4}\\
=& \nabla_x\tilde u_0\nabla_x\tilde \mu_1+\nabla_x\tilde u_1\nabla_x\tilde\mu_0+(\hat\chi_0\tilde u_{0z}\tilde \mu_{3z}+(\hat \chi_1\tilde u_{0z}\varepsilon^{-1}+\tilde u_{1z})\tilde\mu_{2z}\notag\\
&+(\hat \chi_2\tilde u_{0z}\varepsilon^{-2}+\tilde u_{2z})\tilde\mu_{1z}+(\hat \chi_3\tilde u_{0z}\varepsilon^{-3}+\tilde u_{3z})\tilde\mu_{0z})\notag\\
&+\chi_0(1+\tilde u_0)(\Delta_x\tilde\mu_1+\kappa\tilde\mu_{2z}+\tilde\mu_{3zz})\notag\\
&+((1+\tilde u_0)\chi_1\varepsilon^{-1}+\tilde u_1)(\Delta_x\tilde\mu_0+\kappa\tilde\mu_{0z}+\tilde \mu_{2zz})\notag\\
&+((1+\tilde u_0)\chi_2\varepsilon^{-2}+\tilde u_2)(\kappa\tilde \mu_{0z}+\tilde \mu_{1zz})+((1+\tilde u_0)\chi_3\varepsilon^{-3}+\tilde u_3)\tilde\mu_{0zz},\notag
\end{align}
where the solutions of the first and the second equations, following the same treatment as in former sections, derive
\begin{equation}
\tilde u_0 = \tanh(z/\sqrt2),\quad \tilde \mu_0 = 0,\quad \tilde \mu_1 = \kappa S/2,\quad \tilde u_1 = \kappa \Phi.\label{results}
\end{equation}
As for $\tilde \mu_2$, we simplify the equation \eqref{eq-u3-01} by \eqref{results} to derive
\begin{align}
0=\hat \chi_0 \tilde u_{0z}\mu_{2z}+\chi_0 (1+\tilde u_0)\tilde \mu_{2zz},
\end{align}
which leads to $\tilde\mu_2 = b_2$ in $(-\eta,+\infty)$, where $b_2$ is a constant independent of $z$. Recall that by asymptotic matching,
\begin{equation}
\mu_2 = F''(\tilde u_0)\tilde u_2 - \tilde u_{2zz} - \kappa\tilde u_{1z}+F'''(\tilde u_0)\tilde u_1^2/2.\label{eq-u3-02}
\end{equation}
Multiplying the equation \eqref{eq-u3-02} by $\tilde u_0'$ and integrating the resulting one over $(-\infty,+\infty)$, we get
\begin{equation}
\tilde \mu_2\tilde u_{0}|_{-\infty}^{+\infty} = -\kappa^2\int_{-\infty}^{+\infty}(\Phi'-\frac12F'''(\tilde u_0)\Phi^2)u_{0z}dz
\end{equation}
in $(\hat \eta,+\infty)$, which is
\begin{equation}\tilde\mu_2 = -\kappa^2S_1/2\label{eq-u3-03}\end{equation} if we let $S_1 = \kappa^2\int_{-\infty}^{+\infty}(\Phi'-\frac12F'''(\tilde u_0)\Phi^2)u_{0z}dz$. Then, as in \cite{dai2012motion}, we extrapolate a little bit and one may assume that $\tilde\mu_2 = -\kappa^2S_1/2$ in $(-\eta,\hat\eta)$.

Now we solve for $\partial_t^\alpha\phi$. Using (\ref{results}), we have
\begin{align}
&\frac{1}{\Gamma(1-\alpha)}\int_0^{t_2}\frac{\partial_{\tau}\phi(x,\tau)\partial_z \tilde u_0(\phi(x,\tau)/\varepsilon)\chi_0}{(t_2-\tau)^\alpha}d\tau\\
= &\tilde u_{1z}\tilde \mu_{2z}+\tilde u_1\tilde \mu_{2zz}\chi_0(1+\tilde u_0)\tilde\mu_{3zz}+\hat\chi_0\tilde u_{0z}\tilde \mu_{3z}\notag\\
= &\chi_0\partial_z((1+\tilde u_0)\tilde\mu_{3z})-(\chi_0-\hat\chi_0)\tilde u_{0z}\tilde\mu_{3z}+\partial_z(\tilde u_{1}\tilde \mu_{2z}),\notag
\end{align} 
which yields, by integrating over $(-\infty,+\infty)$, that
\begin{align}
\partial_{t_2}^\alpha\phi(2+O(\varepsilon)) =& \lim_{-\eta\rightarrow-\infty}((1+\tilde u_0)\tilde\mu_{3z})+2\partial_\textbf{m}\mu_2^++\lim_{-\eta\rightarrow-\infty}(\tilde u_{1}\tilde \mu_{2z})-\lim_{\hat\eta\rightarrow+\infty}u_{0z}\tilde\mu_{3z}+O(\varepsilon).
\end{align}
Using the matching conditions and letting $\varepsilon\rightarrow 0$, we get
\begin{align}
2\partial_{t_2}^\alpha\phi=2\partial_\textbf{m}\mu_2^++u_1^-\partial_\textbf{m}\mu_1^-,
\end{align}
which gives, by using $\mu_1 = F''(u_0)u_1 = 2u_1$, that
\begin{align}
\partial_{t_2}^\alpha\phi=\partial_\textbf{m}\mu_2^++\frac14\mu_1^-\partial_\textbf{m}\mu_1^-.\label{eq-u3-04}
\end{align}
\textbf{Sharp interface model at $t_2 = \varepsilon^{\frac{2}{\alpha}}t$.}
Combining \eqref{outer31} and \eqref{outer32} with \eqref{results},\eqref{eq-u3-03} and \eqref{eq-u3-04}, we finally derive the sharp interface model in this timescale
\begin{subequations}
\begin{align}
&& &\nabla(\mu_1\nabla \mu_1) = \partial_{t_2}^\alpha u_0,& &\text{in} \quad \Omega^-,&\\
&& &\mu_1 = -\kappa \frac{S}{[U]},& &\text{on}\quad \Gamma,&\\
&& &2\Delta \mu_2 = \partial_{t_2}^\alpha u_0,& &\text{in} \quad \Omega^+,&\label{sharpmullin22}\\
&& &\mu_2 = -\kappa^2\frac{S_1}{[U]},& &\text{on}\quad \Gamma,&\\
&& &\text{I}^{1-\alpha}V =\partial_\textbf{m}\mu_2^++\frac14\mu_1^-\partial_\textbf{m}\mu_1^-,& &\text{on}\quad \Gamma.&
\end{align}
\end{subequations}
$u_0$ is the sign function of $\phi$, i.e., $u_0 = \pm1$ in $\Omega^\pm$.


\section{Scaling invariant property and coarsening rate heuristic}\label{sec4} 
In physics, coarsening is a progress when the pattern formed by the material ``coarsens" and during which the ``typical length scale" of the system is increasing. For phase field models, since the energy of the system is proportional to the area of the interfacial layer, energy decay would result in the reduction of the interface layer and the pattern coarsens. Coarsening phenomena are also observed in numerical simulations of the pattern formation governed by TFCHE. As many people believe, coarsening is due to some ``scaling invariant" property of the system, so the scaling-invariant power law of sharp interface model coincides in the coarsening rate in the simulation of \cite{dai2012motion,kohn2004coarsening,pego2007lectures}.

Considering the nonlocal MS model with constant mobility in $t_1 = \varepsilon^{\frac{1}{\alpha}}t$ time scale. It is scaling invariant in the following sense. Rescaling $\mu$, $x$, $t$ and $\phi$ by
\begin{displaymath}
 x = X\hat{x}, \quad  t = T\hat{t},\quad \mu = M\hat{\mu}_1, \quad \phi(x,t) = X\hat{\phi}(\hat{x},\hat{t}).
\end{displaymath}
Direct calculation leads to
\begin{displaymath}
\kappa = X^{-1}\hat\kappa,\quad \partial^\alpha_t\phi = X/T^\alpha\partial^\alpha_{\hat t}\hat\phi, \quad \partial_\textbf{m}\mu =X^{-2} \partial_{\hat{m}}\hat\mu,
\end{displaymath} 
and
\begin{subequations}
\begin{align}
&& &\frac{M}{X^2}\hat\Delta \hat\mu = \frac{1}{T^\alpha}\partial^\alpha_{\hat t} u_0,  & &\text{in} \quad \Omega\backslash \Gamma,&\label{coarsening-eq-1-start}\\
&& &M\hat\mu = \frac{1}{X}\hat\kappa \frac{S}{[U]},& &\text{on}\quad \Gamma,&\\
&& &\frac{X}{T^\alpha}\partial^\alpha_{\hat t}\hat\phi = \frac{M}{X}[ \partial_{\hat{m}}\hat\mu]_-^+[U]^{-1} , & &\text{on}\quad \Gamma&\label{coarsening-eq-1-end}.
\end{align}
\end{subequations}
If taking $M = X^{-1}$ and $T^\alpha = X^3$, the system has exactly the same form as (\ref{sharpmullin21})-(\ref{sharpmullin24}). This is the scaling invariance property and it shows that the typical length scale $l$ of this model satisfies a $l\sim ct^{\frac\alpha3}$ power law, which implies that the TFCHE admits a coarsening rate of $\frac\alpha3$. This result fits the numerical experiments in \cite{tang2019energy,zhao2019power} well. 

In the second part of this section, we aim to use this idea to determine the coarsening rate of the sharp interface models of the degenerate TFCHE.
Firstly, for the sharp interface models in $t_1 = \varepsilon^{\frac{1}{\alpha}}t$ time scale
\begin{subequations}
\begin{align}
 && & \partial_{t_1}^\alpha u_0 =\Delta \mu_1,& & \text{in} \quad\Omega^+,&\label{stefan-1+u-start}\\
&& & \mu_1 = -\kappa \frac{S}{[U]},&  &\text{on}\quad \Gamma,&\\
&& & \text{I}^{1-\alpha}V = \partial_\textbf{m}\mu_1^+,& &\text{on}\quad \Gamma,&\label{stefan-1+u-end}
\end{align} 
\end{subequations}
by using the same rescaling as in \eqref{coarsening-eq-1-end} - \eqref{coarsening-eq-1-start}
\begin{displaymath}
 x = \lambda^{\frac\alpha3}\hat{x}, \quad  t_1 = \lambda\hat{t}_1,\quad \mu_1 = \lambda^{-\frac\alpha3}\hat{\mu}_1, \quad \phi(x,t_1) =  \lambda^{\frac\alpha3}\hat{\phi}(\hat{x},\hat{t}_1),
\end{displaymath}
we find that \eqref{stefan-1+u-start}-\eqref{stefan-1+u-end} preserves a ${\frac\alpha3}$ coarsening rate, too.

On the other hand, for the sharp interface model \eqref{dgms21}-\eqref{dgms22} in $t_2 = \varepsilon^{\frac{2}{\alpha}}t$ time scale
\begin{subequations}
\begin{align}
 &&  & \partial_{t_2}^\alpha u_0  = \nabla(\mu_1\nabla \mu_1),& & \text{in} \quad \Omega^-,&\label{eq-si-1}\\
&& &\mu_1 = -\kappa \frac{S}{[U]},& &\text{on}\quad \Gamma,&\\
&& &\partial_{t_2}^\alpha u_0  = 2\Delta \mu_2,& &\text{in} \quad \Omega^+,&\\
&& &\mu_2 = -\kappa^2\frac{S_1}{[U]},& & \text{on}\quad \Gamma,&\\
&& &\text{I}^{1-\alpha}V =\partial_\textbf{m}\mu_2^++\frac14\mu_1^-\partial_\textbf{m}\mu_1^-,& &\text{on}\quad \Gamma\label{eq-si-end}.
\end{align} 
\end{subequations}
Taking $M_1, M_2, T, X$ to be the length scales of the chemical potentials, time, and space, respectively, we rescale the above system \eqref{eq-si-1}-\eqref{eq-si-end} so that
\begin{subequations}
\begin{align}
&& &\frac{1}{T^\alpha}\partial_{\hat t}^\alpha u_0  = \frac{M_1^2}{X^2}\nabla(\hat\mu_1\nabla \hat\mu_1),& &\text{in} \quad \Omega^-,&\\
&& &M_1\hat\mu_1 = -\frac{1}{X}\hat\kappa \frac{S}{[U]},& &\text{on}\quad \Gamma,&\\
&& &\frac{1}{T^\alpha}\partial_{\hat t}^\alpha u_0 = 2\frac{M_2}{X^2}\hat\Delta \hat\mu_2,& &\text{in} \quad \Omega^+,&\\
&& &M_2\hat\mu_2 = -\frac{1}{X^2}\hat\kappa^2\frac{S_1}{[U]},& &\text{on}\quad \Gamma,&\\
&& &\frac{X}{T^\alpha}\partial^\alpha_{\hat t}\text{I}^{1-\alpha}\hat V =\frac{M_2}{X}\partial_\textbf{m}\hat\mu_2^++\frac{M_1^2}{X}\frac14\hat\mu_1^-\partial_\textbf{m}\hat\mu_1^-,& &\text{on}\quad \Gamma.&
\end{align} 
\end{subequations}
The system is the same form as \eqref{dgms21}-\eqref{dgms22} if we take
\begin{displaymath}
 T^\alpha = X^4, \quad  M_1 = \frac{1}{X},\quad \mathrm{and}\quad M_2 =\frac{1}{X^2}.
\end{displaymath}
It exhibits a power law relation $l\sim ct^{\frac{\alpha}{4}}$. Moreover, this power law indicate a coarsening rate of $\frac\alpha4$.

\section{Conclusions}
We study the front motion and obtain the corresponding sharp interface models of the TFCHE with two different kinds of diffusion mobilities. We find that in both cases the sharp interface limits are sensitive to the timescale. For example, in a slow time scale $\varepsilon^{\frac1\alpha}t$, the asymptotic limits are fractional Mullins--Sekerka(MS) models, which are formally similar to classical MS models excepted for the non-local term. 

Moreover, power-law arguments show that the nonlocal fractional MS model of TFCHE with constant mobility fits the $\frac\alpha3$ coarsening rate obtained in existing numerical experiments \cite{zhao2019power,tang2019energy}. Moreover, TFCHE with the one-sided degenerate mobility contains two stages of different coarsening rates $\frac{\alpha}{3}$ and $\frac{\alpha}{4}$. The results show that the TFCHE might could be use to model the coarsening process with a general coarsening rate. We expect to extend similar arguments to the nonlocal-in-time phase-field equations, in which the time fractional operator is replaced by a nonlocal-in-time operator \cite{du2017analysis}.

\subsection*{Acknowledgement}
This work is partially supported by the Special Project on High-Performance Computing of the National Key R\&D Program under No. 2016YFB0200604, the National Natural Science Foundation of China (NSFC) Grant No. 11731006, the NSFC/Hong Kong RGC Joint Research Scheme (NSFC/RGC 11961160718), and the fund of the Guangdong Provincial Key Laboratory of Computational Science and Material Design (No. 2019B030301001). The work of J. Yang is supported by the National Science Foundation of China (NSFC-11871264).

\bibliographystyle{siam}
\bibliography{refers}
\end{document}